\documentclass[12pt]{article}
\usepackage{latexsym,amsfonts,amssymb}
\usepackage[dvips]{color}
\usepackage{colordvi,multicol}
\usepackage{amsmath}

\newtheorem{thm}{Theorem}[section]

\newtheorem{pro}[thm]{Proposition}

\newtheorem{rem}[thm]{Remark}

\date{}
\begin{document}

\title{\bf An Opposite Gaussian Product Inequality}
 \author{Oliver Russell and Wei Sun\\ \\ \\
  {\small Department of Mathematics and Statistics}\\
    {\small Concordia University, Canada}\\ \\
{\small o\_russel@live.concordia.ca,\ \ \ \ wei.sun@concordia.ca}}

\maketitle

\begin{abstract}

\noindent The long-standing Gaussian product inequality (GPI) conjecture states that
$E [\prod_{j=1}^{n}|X_j|^{\alpha_j}]\geq\prod_{j=1}^{n}E[|X_j|^{\alpha_j}]$ for any centered Gaussian random vector $(X_1,\dots,X_n)$ and any non-negative real numbers $\alpha_j$, $j=1,\ldots,{n}$. In this note, we prove a novel ``opposite GPI" for centered bivariate Gaussian random variables when $-1<\alpha_1<0$ and $\alpha_2>0$: $E[|X_1|^{\alpha_1}|X_2|^{\alpha_2}]\le E[|X_1|^{\alpha_1}]E[|X_2|^{\alpha_2}]$. This completes the picture of bivariate Gaussian product relations.
\end{abstract}



\section{Introduction}

Multivariate Gaussian distributions are essential to the theory and applications of probability, statistics and many other branches of science. Since Royen's 2014 proof of the long-conjectured Gaussian correlation inequality \cite{Royen,LM}, interest in inequalities related to multivariate Gaussian distributions has risen drastically in the research community. In particular, the focus has shifted towards solving the long-standing Gaussian product inequality (GPI) conjecture, which was stated most generally by Li and Wei \cite{LW12} as follows: for any non-negative real numbers $\alpha_j$, $j=1,\ldots,{n}$, and  any ${n}$-dimensional real-valued centered Gaussian random vector $(X_1,\dots,X_{n})$,
\begin{eqnarray}\label{LW-inequ}
E \left[\prod_{j=1}^{n}|X_j|^{\alpha_j}\right]\geq \prod_{j=1}^{n}E[|X_j|^{\alpha_j}].
\end{eqnarray}

The verification of the GPI would have immediate and significant impact on multiple problems in other fields. For example (see Malicet et al. \cite{MNPP16}), if (\ref{LW-inequ}) holds when the $\alpha_j$'s are any equal positive  even integers, then the `real linear polarization constant' conjecture raised by Ben\'{\i}tem  et al. \cite{BST98} is true, and (\ref{LW-inequ}) is deeply linked to the celebrated $U$-conjecture by Kagan, Linnik and Rao \cite{KLR1973}.

Due to the significant difficulty of proving the GPI conjecture, only partial results have been obtained so far. In \cite{Fr08}, Frenkel used algebraic methods to prove (\ref{LW-inequ}) for the case that $\alpha_j=2$ for each $j=1,\ldots,{n}$. In \cite[Theorem 3.2]{LHS}, Lan et al. used the Gaussian hypergeometric functions to prove the following GPI: for any $m_1,m_2\in\mathbb{N}$ and any centered Gaussian random vector $(X_1,X_2,X_3)$,
$$
E[X_1^{2m_1}X_2^{2m_2}X_3^{2m_2}]\ge E[X_1^{2m_1}]E[X_2^{2m_2}]E[X_3^{2m_2}].
$$
Most recently, when the $\alpha_j$'s are positive even integers, we have proved more $3$-dimensional GPIs \cite{RusSun}, designed a computational algorithm based on sums-of-squares polynomials to verify any specific case of the GPI \cite{RusSun2}, and completely solved the GPI when all correlations are non-negative \cite[Lemma 2.3]{RusSun}. Genest and Ouimet \cite{GenOui} gave a combinatorial proof of the GPI under the stronger assumption that there exists a matrix $C\in[0,\infty)^{n\times n}$ such that $(X_1,\dots,X_n)=(U_1,\dots,U_n)C$ in law, where $(U_1,\dots,U_n)$
is an $n$-dimensional standard Gaussian random vector. Additionally, Edelmann et al. \cite{Edel} used a different method to extend \cite[Lemma 2.3]{RusSun} to the multivariate gamma distribution.

Despite all of this progress, solving the GPI conjecture is still very challenging when correlations can be negative, or when the $\alpha_j$'s are not necessarily positive even integers. However, for the $2$-dimensional case, there are more general results. In \cite{We14}, Wei used
integral representations to prove a stronger version of (\ref{LW-inequ}) for $\alpha_j\in (-1,0)$ as follows:
\begin{eqnarray*}
E \left[\prod_{j=1}^{n}|X_j|^{\alpha_j}\right]\geq
E \left[\prod_{j=1}^k|X_j|^{\alpha_j}\right]E\left[\prod_{j=k+1}^{n}|X_j|^{\alpha_j}\right],\ \ \ \ \forall 1\le k\le n-1.
\end{eqnarray*}
By Karlin and Rinott \cite[Corollary 1.1, Theorem 3.1 and Remark 1.4]{KR}, (\ref{LW-inequ}) holds for $n=2$ by virtue of $(|X_1|,|X_2|)$ possessing the MTP$_2$ property. Thus, we know the following $2$-dimensional GPI holds:
\begin{thm}\label{jkl} Let $(X_1,X_2)$ be centered bivariate Gaussian random variables and $\alpha_1,\alpha_2>0$ or $-1<\alpha_1,\alpha_2<0$.  Then,
$$
E[|X_1|^{\alpha_1}|X_2|^{\alpha_2}]\ge E[|X_1|^{\alpha_1}]E[|X_2|^{\alpha_2}],
$$
and the equality sign holds if and only if $X_1,X_2$ are independent.
\end{thm}

In this note, we complete the picture of $2$-dimensional Gaussian product relations by proving what we call an ``opposite GPI" when $\alpha_1$ and $\alpha_2$ are of opposite sign. In Section 2, using the hypergeometric function, we state and prove this main result (see Theorem \ref{thm2}) and then give a more explicit expression for the $2$-dimensional GPI and opposite GPI when one of $\alpha_1,\alpha_2$ is a positive even integer. In section 3, we prove a $1$-dimensional GPI and opposite GPI. The results contained in this note extend and conclude the study of the $2$-dimensional GPI in its full generality.

\section{Opposite GPI}\setcounter{equation}{0}

Throughout this note, any Gaussian random variable is assumed to be real-valued and non-degenerate, i.e., has positive variance. Our main result is the following novel opposite GPI:

\begin{thm}\label{thm2} Let $(X_1,X_2)$ be centered bivariate Gaussian random variables, $-1<\alpha_1<0$ and $\alpha_2>0$.  Then,
$$
E[|X_1|^{\alpha_1}|X_2|^{\alpha_2}]\le E[|X_1|^{\alpha_1}]E[|X_2|^{\alpha_2}],
$$
and the equality sign holds if and only if $X_1,X_2$ are independent.
\end{thm}

\noindent Proof.\ \  Let $U$ be a standard Gaussian random variable and $\nu>-1$. We have
\begin{eqnarray}\label{dfg}
E[|U|^{\nu}]&=&\frac{\sqrt{2}}{\sqrt{\pi}}\int_0^{\infty}x^{\nu}e^{-\frac{x^2}{2}}dx\nonumber\\
&=&\frac{2^{\nu/2}}{\sqrt{\pi}}\int_0^{\infty}z^{(\nu-1)/2}e^{-z}dz\nonumber\\
&=&\frac{2^{\nu/2}\Gamma(\frac{\nu+1}{2})}{\sqrt{\pi}}.
\end{eqnarray}

We assume without loss of generality that  ${\rm Var}(X_1)={\rm Var}(X_2)=1$. Denote by $\rho$ the correlation of $X_1$ and $X_2$. By (\ref{dfg}), we get
\begin{eqnarray}\label{gh2}
E[|X_1|^{\alpha_1}]E[|X_2|^{\alpha_2}]=\frac{2^{(\alpha_1+\alpha_2)/2}\Gamma(\frac{\alpha_1+1}{2})\Gamma(\frac{\alpha_2+1}{2})}{{\pi}}.
\end{eqnarray}
If $|\rho|<1$, by Nabeya \cite{N51}, we have that
\begin{eqnarray}\label{gh3}
E[|X_1|^{\alpha_1}|X_2|^{\alpha_2}]=\frac{2^{(\alpha_1+\alpha_2)/2}\Gamma(\frac{\alpha_1+1}{2})\Gamma(\frac{\alpha_2+1}{2})}{\pi}\cdot F\left(-\frac{\alpha_1}{2},-\frac{\alpha_2}{2};\frac{1}{2};\rho^2\right),
\end{eqnarray}
where $F$ is a Gaussian hypergeometric function given by
$$
F(\alpha,\beta;\gamma;z)=1+ \frac{\alpha\beta}{1!\gamma}z+\frac{\alpha(\alpha+1)
\beta(\beta+1)}{2!\gamma(\gamma+1)}z^2+\cdots,\ \ \ \ |z|<1.
$$

Define
$$
G(z):= F\left(-\frac{\alpha_1}{2},-\frac{\alpha_2}{2};\frac{1}{2};z\right),\ \ \ \ |z|<1.
$$
By  Bateman \cite[2.8-(20), page 102]{B} and the Euler transformation (cf. Rainville \cite[Chapter 4, Theorem 21, page 60]{R}), we get
\begin{eqnarray}\label{vbn}
G'(z)&=&\frac{\alpha_1\alpha_2}{2}\cdot F\left(1-\frac{\alpha_1}{2},1-\frac{\alpha_2}{2};\frac{3}{2};z\right)\nonumber\\
&=&\frac{\alpha_1\alpha_2}{2}\cdot(1-z)^{(\alpha_1+\alpha_2-1)/2} F\left(\frac{\alpha_1+1}{2},\frac{\alpha_2+1}{2};\frac{3}{2};z\right)\nonumber\\
&<&0,\ \ \ \ |z|<1.
\end{eqnarray}
Therefore, the proof is complete by (\ref{gh2}) and (\ref{gh3}). \hfill \fbox

\begin{rem}
Following the proof of Theorem \ref{thm2}, we can give a new proof of Theorem \ref{jkl} using the hypergeometric function. The main difference is that now $\alpha_1\alpha_2>0$ and hence (\ref{vbn}) is replaced by the inequality $G'(z)>0$, $|z|<1$.
\end{rem}

If the exponent $\alpha_2$ is a positive even integer, then we can get a more explicit comparison of the joint moment and the product of marginal moments. Obviously, this comparison result gives a different proof for the GPI (when $\alpha_1>0$) and the opposite GPI (when $-1<\alpha_1<0$).

\begin{pro}
Let $(X_1,X_2)$ be centered bivariate Gaussian random variables with correlation $\rho$,  $\alpha_1>-1$ and $m\in\mathbb{N}$. Then,
\begin{eqnarray*}
&&\frac{E[|X_1|^{\alpha_1}|X_2|^{2m}]}{E[|X_1|^{\alpha_1}]E[|X_2|^{2m}]}\\
&=&(1-\rho^2)^m+\sum_{j=1}^m{{{m}\choose j}}{\rho^{2j}}{(1-\rho^2)}^{m-j}\frac{[\alpha_1+(2j-1)][\alpha_1+(2j-3)]\cdots(\alpha_1+1)}{(2j-1)!!}.
\end{eqnarray*}
\end{pro}

\noindent Proof.\ \ We assume without loss of generality that $X_1=U_1$ and $X_2=aU_1+U_2$, where $U_1,U_2$ are independent  standard Gaussian random variables and $a\in\mathbb{R}$.

By (\ref{dfg}), we get
\begin{eqnarray}\label{zxcv}
E[|X_1|^{\alpha_1}]E[|X_2|^{2m}]=E[|U_1|^{\alpha_1}]E[|aU_1+U_2|^{2m}]=\frac{2^{\alpha_1/2}\Gamma\left(\frac{\alpha_1+1}{2}\right)(1+a^2)^m(2m-1)!!}{\sqrt{\pi}},
\end{eqnarray}
and
\begin{eqnarray}\label{zxcv2}
&&E[|X_1|^{\alpha_1}|X_2|^{2m}]\nonumber\\
&=&E[|U_1|^{\alpha_1}|aU_1+U_2|^{2m}]\nonumber\\
&=&\frac{1}{2\pi}\int_{-\infty}^{\infty}\int_{-\infty}^{\infty}|x_1|^{\alpha_1}|ax_1+x_2|^{2m}e^{-\frac{x_1^2+x_2^2}{2}}dx_1dx_2\nonumber\\
&=&\frac{1}{2\pi}\int_{-\infty}^{\infty}|x_1|^{\alpha_1}e^{-\frac{x_1^2}{2}}\int_{-\infty}^{\infty}\left[\sum_{i=0}^{2m}{{2m}\choose i}(ax_1)^ix_2^{2m-i}\right]e^{-\frac{x_2^2}{2}}dx_2dx_1\nonumber\\
&=&\sum_{j=0}^{m-1}{{2m}\choose 2j}(2m-2j-1)!!\frac{1}{\sqrt{2\pi}}\int_{-\infty}^{\infty}a^{2j}|x_1|^{\alpha_1+2j}e^{-\frac{x_1^2}{2}}dx_1\nonumber\\
&&+\frac{1}{\sqrt{2\pi}}\int_{-\infty}^{\infty}a^{2m}|x_1|^{\alpha_1+2m}e^{-\frac{x_1^2}{2}}dx_1\nonumber\\
&=&(2m-1)!!\left\{\frac{1}{\sqrt{2\pi}}\int_{-\infty}^{\infty}|x_1|^{\alpha_1}e^{-\frac{x_1^2}{2}}dx_1+\sum_{j=1}^{m}\frac{{{m}\choose j}a^{2j}}{(2j-1)!!}\cdot\frac{1}{\sqrt{2\pi}}\int_{-\infty}^{\infty}|x_1|^{\alpha_1+2j}e^{-\frac{x_1^2}{2}}dx_1\right\}\nonumber\\
&=&\frac{(2m-1)!!}{\sqrt{\pi}}\left\{2^{\alpha_1/2}\Gamma\left(\frac{\alpha_1+1}{2}\right)+\sum_{j=1}^{m}\frac{{{m}\choose j}a^{2j}2^{(\alpha_1+2j)/2}\Gamma(\frac{\alpha_1+2j+1}{2})}{(2j-1)!!}\right\}\nonumber\\
&=&\frac{2^{\alpha_1/2}\Gamma\left(\frac{\alpha_1+1}{2}\right)(2m-1)!!}{\sqrt{\pi}}\left\{1+\sum_{j=1}^m{{{m}\choose j}a^{2j}}\frac{[\alpha_1+(2j-1)][\alpha_1+(2j-3)]\cdots(\alpha_1+1)}{(2j-1)!!}\right\}.\ \ \ \ \ \ \ \ \ \ 
\end{eqnarray}
Note that
$$
a^2=\frac{\rho^2}{1-\rho^2},\ \ \ \ \ \ 1+a^2=\frac{1}{1-\rho^2}.
$$
Therefore, the proof is complete by (\ref{zxcv}) and (\ref{zxcv2}).\hfill\fbox

\section{1-Dimensional Inequalities}\setcounter{equation}{0}

In this section, we present a $1$-dimensional GPI and opposite GPI.

\begin{pro}
Let $X$ be a centered Gaussian random variable.

\noindent (i) If $-1<\alpha_1<0$ and $\alpha_2>0$, then
$$
E[|X|^{\alpha_1+\alpha_2}]<\frac{(\alpha_1+1)(\alpha_2+1)}{\alpha_1+\alpha_2+1}E[|X|^{\alpha_1}]E[|X|^{\alpha_2}].
$$

\noindent (ii) If $\alpha_1,\alpha_2>0$ or $-1<\alpha_1,\alpha_2<0$ with $\alpha_1+\alpha_2>-1$, then
$$
E[|X|^{\alpha_1+\alpha_2}]>\frac{(\alpha_1+1)(\alpha_2+1)}{\alpha_1+\alpha_2+1}E[|X|^{\alpha_1}]E[|X|^{\alpha_2}].
$$
\end{pro}

\noindent Proof.\ \  We assume without loss of generality that  ${\rm Var}(X)=1$. By (\ref{dfg}), we get
$$
E[|X|^{\alpha_1+\alpha_2}]=\frac{2^{(\alpha_1+\alpha_2)/2}\Gamma(\frac{\alpha_1+\alpha_2+1}{2})}{\sqrt{\pi}},
$$
and
\begin{eqnarray*}
E[|X|^{\alpha_1}]E[|X|^{\alpha_2}]=\frac{2^{(\alpha_1+\alpha_2)/2}\Gamma(\frac{\alpha_1+1}{2})\Gamma(\frac{\alpha_2+1}{2})}{{\pi}}.
\end{eqnarray*}
Hence,
\begin{eqnarray}\label{as}
\frac{E[|X|^{\alpha_1+\alpha_2}]}{E[|X|^{\alpha_1}]E[|X|^{\alpha_2}]}
=\frac{\Gamma(\frac{\alpha_1+\alpha_2+1}{2})\Gamma(\frac{1}{2})}{\Gamma(\frac{\alpha_1+1}{2})\Gamma(\frac{\alpha_2+1}{2})}
=\frac{B(\frac{\alpha_1+\alpha_2+1}{2},\frac{1}{2})}{B(\frac{\alpha_1+1}{2},\frac{\alpha_2+1}{2})}.
\end{eqnarray}

Note that (cf. \cite[p. 8]{AAR})
$$
B(x,y)=\frac{x+y}{xy}\prod_{n=1}^{\infty}\left(1+\frac{xy}{n(x+y+n)}\right)^{-1},\ \ \ \ x,y>0,
$$
and
\begin{eqnarray*}
&&\frac{\alpha_1+\alpha_2+1}{2}\cdot\frac{1}{2}>\frac{\alpha_1+1}{2}\cdot\frac{\alpha_2+1}{2},\ \ \ \ {\rm if}\ -1<\alpha_1<0\ {\rm and}\ \alpha_2>0,\\
&&\frac{\alpha_1+\alpha_2+1}{2}\cdot\frac{1}{2}<\frac{\alpha_1+1}{2}\cdot\frac{\alpha_2+1}{2},\ \ \ \ {\rm if}\ \alpha_1,\alpha_2>0\ {\rm or}\ -1<\alpha_1,\alpha_2<0.
\end{eqnarray*}
Then,
\begin{eqnarray*}
&&\frac{B(\frac{\alpha_1+\alpha_2+1}{2},\frac{1}{2})}{B(\frac{\alpha_1+1}{2},\frac{\alpha_2+1}{2})}<\frac{(\alpha_1+1)(\alpha_2+1)}{\alpha_1+\alpha_2+1},\ \ \ \ {\rm if}\ -1<\alpha_1<0\ {\rm and}\ \alpha_2>0,\\
&&\frac{B(\frac{\alpha_1+\alpha_2+1}{2},\frac{1}{2})}{B(\frac{\alpha_1+1}{2},\frac{\alpha_2+1}{2})}>\frac{(\alpha_1+1)(\alpha_2+1)}{\alpha_1+\alpha_2+1},\ \ \ \ {\rm if}\ \alpha_1,\alpha_2>0\ {\rm or}\ -1<\alpha_1,\alpha_2<0\ {\rm with}\ \alpha_1+\alpha_2>-1.
\end{eqnarray*}
Therefore, the proof is complete by (\ref{as}).\hfill \fbox

\section{Discussion}

As mentioned in the introduction, the bivariate GPI had been proved for exponents with the same signs, but we could find no prior study of cases with opposite signs. In this note, we used the Gaussian hypergeometric function to prove an opposite GPI, which we believe to be the first of its kind. Furthermore, our method provides an extremely simple alternative proof to Theorem \ref{jkl}. Therefore, not only do our new results fill in the last remaining piece of the $2$-dimensional puzzle, but this note also serves as a self-contained complete picture.
\hfill \break

\begin{large} \noindent\textbf{Acknowledgements} \end{large} This work was supported by the Natural Sciences and Engineering Research Council of Canada (Nos. 559668-2021 and
4394-2018).

\end{document}